\documentclass[a4paper,UKenglish,cleveref, autoref]{lipics-v2019}

\usepackage{graphicx} 
\usepackage{verbatim} 
\usepackage{amsmath,amssymb,amsthm,bbm}
\usepackage{color}
\usepackage{enumerate}
\usepackage{array}
\usepackage[noadjust]{cite}
\usepackage{multirow}

\usepackage{xcolor}
\usepackage{hyperref}
\definecolor{darkgreen}{rgb}{0,0.4,0}
\definecolor{BrickRed}{rgb}{0.65,0.08,0}
\hypersetup{colorlinks=true,linkcolor=blue,citecolor=red,filecolor=BrickRed,urlcolor=darkgreen}
\nolinenumbers

\DeclareMathOperator\Pol{Pol}

\newcommand{\LandauO}{\mathcal{O}}

\newcommand{\Cc}{\mathcal{C}}

\newcommand{\Qc}{\mathcal{Q}}

\newcommand{\Sc}{\mathcal{S}}

\newcommand{\C}{\mathbb{C}}

\newcommand{\Q}{\mathbb{Q}}

\newcommand{\R}{{\mathbb R}}

\newcommand{\ns}{{\mathbb N}} 
\newcommand{\zs}{{\mathbb Z}} 
\newcommand{\qs}{{\mathbb Q}}  

\newcommand{\beq}{\begin{equation}}
\newcommand{\eeq}{\end{equation}}

\def\emm#1,{{\em #1}}

\newcommand{\GA}{\mathbb{A}}
\newcommand{\cS}{\mathcal S}

\newcommand{\bx}{\bar{x}}
\newcommand{\by}{\bar{y}}

\newcommand{\uroot}{\zeta}  
\newcommand{\spol}{S}

\newcommand{\xzero}{u}
\newcommand{\av}{v}
\newcommand{\aw}{w}

\newcommand{\DSA}{B_1}
\newcommand{\DSB}{B_2}

\newcommand{\hookSEarrow}{\mathrel{\rotatebox[origin=c]{-45}{$\hookrightarrow$}}}
\newcommand{\hookNEarrow}{\mathrel{\rotatebox[origin=c]{45}{$\hookrightarrow$}}}

\newcommand{\gf}{generating function}
\newcommand{\gfs}{generating functions}
\newcommand{\fps}{formal power series}

\title{More models of walks avoiding a quadrant (extended abstract)} 

\titlerunning{More models of walks avoiding a quadrant} 

\author{Mireille Bousquet-M{\'e}lou}{CNRS, Universit\'e de Bordeaux, Laboratoire Bordelais de Recherche en Informatique, UMR 5800, 351
  cours de la Libération, 33405 Talence Cedex, France}{mireille.bousquet-melou@u-bordeaux.fr}{}{}

\author{Michael Wallner}{Universit\'e de Bordeaux, Laboratoire Bordelais de Recherche en Informatique, UMR 5800, 351 
  cours de la Libération, 33405 Talence Cedex, France \and TU Wien, Institute for Discrete Mathematics and Geometry, Wiedner Hauptstra{\ss}e 8--10, 1040 Wien, Austria \and \url{https://dmg.tuwien.ac.at/mwallner/}}{michael.wallner@tuwien.ac.at}{https://orcid.org/0000-0001-8581-449X}{%
Supported by the Erwin Schr{\"o}dinger Fellowship of the Austrian Science Fund (FWF):~J~4162-N35.
}

\authorrunning{M. Bousquet-M{\'e}lou and M. Wallner}

\Copyright{Mireille Bousquet-M{\'e}lou and Michael Wallner}



\category{}

\supplement{}

\acknowledgements{We thank our referees for their careful reading.}

\EventEditors{Michael Drmota and Clemens Heuberger}
\EventNoEds{2}
\EventLongTitle{31st International Conference on Probabilistic, Combinatorial and Asymptotic Methods for the Analysis of Algorithms (AofA 2020)}
\EventShortTitle{AofA 2020}
\EventAcronym{AofA}
\EventYear{2020}
\EventDate{June 15--19, 2020}
\EventLocation{Klagenfurt, Austria}
\EventLogo{}
\SeriesVolume{159}
\ArticleNo{13}

\begin{document}

\maketitle

\begin{abstract}
  We continue the enumeration of plane lattice paths avoiding the negative quadrant initiated by the first author in~\cite{Bousquet2016}. We solve in detail a new case, {the king walks,}
  where all $8$ nearest neighbour steps are allowed. As in the two cases solved in~\cite{Bousquet2016}, the associated \gf\ is proved to differ from a simple, explicit D-finite series (related to the enumeration of walks confined to the first quadrant) by an algebraic one. The principle of the approach is the same as in~\cite{Bousquet2016}, but challenging theoretical and computational difficulties arise as we now handle algebraic series of larger degree.

  We also explain why we expect the observed algebraicity phenomenon to persist for $4$ more models, for which the quadrant problem is solvable using the reflection principle.
\end{abstract}

\clearpage


\section{Introduction}
\label{sec:intro}

In this paper we continue the enumeration of plane lattice paths
confined to non-convex cones initiated by the first author in~\cite{Bousquet2016}. 
Therein the two most natural models of walks confined to the three-quadrant cone $\Cc := \{ (i,j) : i \geq 0 \text{ or } j \geq
0 \}$ were studied: walks with steps $\{\rightarrow,  \uparrow,
\leftarrow, \downarrow \}$, and those with steps $\{\nearrow, \nwarrow, \swarrow,  \searrow \}$.
In both cases, the generating function that counts walks starting at the origin
was proved to differ {(additively)}
from a simple explicit {D-finite}
 series by an algebraic one.
The tools essentially involved power series manipulations, coefficient extractions, and polynomial elimination. 

Later, Raschel and Trotignon gave in~\cite{RaschelTrotignon2018Avoiding} {sophisticated}  integral expressions for $8$ models,
which imply that $3$ additional models ($\{\nearrow, \leftarrow,  \downarrow \}$, $\{\rightarrow, \uparrow, \swarrow \}$, and $\{\rightarrow, \nearrow, \uparrow, \leftarrow, \swarrow, \downarrow \}$) are D-finite.
Their results use an analytic approach inspired by earlier work on probabilistic and enumerative aspects of quadrant walks~\cite{fayolle-livre,raschel-unified}.

\begin{figure}[ht]
  \centering
 \scalebox{0.87}{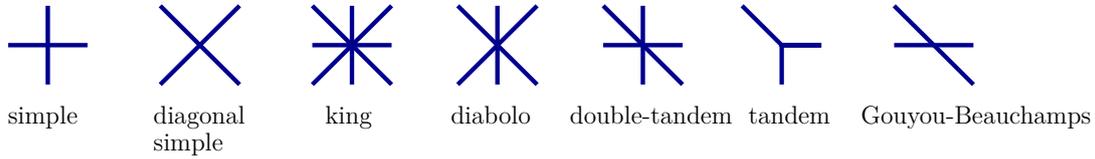}
 \caption{The seven step sets
   to which the strategy of this
    paper should apply. The first two are solved in~\cite{Bousquet2016}, the third one in this paper. }
  \label{fig:7models}
\end{figure}

In this paper we first extend the results of~\cite{Bousquet2016} to
the so-called {\emph{king
walks}}, which take their steps {from}
$\{\rightarrow, \nearrow, \uparrow, \nwarrow, \leftarrow, \swarrow,
\downarrow, \searrow \}$.
{We show that} the \emph{algebraicity phenomenon} of~\cite{Bousquet2016} persists: if
$Q(x,y;t)$ (resp.~$C(x,y;t)$) counts walks starting from the origin that are confined to the non-negative quadrant $\Qc:= \{ (i,j) : i \geq 0 \text{ and } j \geq
0 \}$ (resp.~to the cone $\Cc$)
by the length (variable $t$) and the coordinates of the endpoint
(variables $x,y$), then $C(x,y;t)$ differs from the series
\begin{align*}
	\frac{1}{3}\left( Q(x,y;t) -  Q(1/x,y;t)/x^2 - Q(x,1/y;t)/y^2 \right)
\end{align*}
by an algebraic series, as detailed  in our  main theorem below.
Moreover, we expect a similar property to hold  (with variations on the above
linear combination of the series $Q$)
for the $7$ step sets of Figure~\ref{fig:7models}, related to reflection
groups, and for which the quadrant
problem can be solved using the reflection principle~\cite{gessel-zeilberger}.
However, we also expect the effective solution of these models to be
extremely challenging in computational terms, mostly, because the
relevant algebraic series have very large degree. This is
illustrated by our main theorem below. There, and in the sequel, we use the shorthand
$\bx=1/x$, $\by=1/y$, and omit in the notation the dependencies on~$t$,
writing for instance $Q(x,y)$ instead of $Q(x,y;t)$.

\begin{theorem}\label{theo:king}
  Take the step set $\{ -1,0,1\}^2\setminus\{(0,0)\}$ and let
  $Q(x,y)$ be the \gf\ of lattice walks starting from $(0,0)$ that are confined to the first quadrant $\Qc$ (this series is D-finite and given in~\cite{bomi10}). 
  Then, the \gf\ of  walks starting from $(0,0)$,
confined to $\Cc$, and ending in the first quadrant (resp.~at a negative
abscissa)
is
\beq\label{sol-king}
  \frac 1  3 Q(x,y) + P(x,y), \qquad (\hbox{resp.} -\frac{\bx^2}{3}
    Q(\bx,y) + \bx M(\bx,y)
		),
\eeq
where $P(x,y)$ and $M(x,y)$ are algebraic of degree $216$ over $\qs(x,y,t)$. Of course, the \gf\ of walks ending  at a negative ordinate follows, using the $x/y$-symmetry.

The series $P$ is expressed in terms {of} $M$ by:
 \beq\label{Psol-king}
   P(x,y)=\bx \big( M(x,y)-M(0,y)\big) +\by \big(
     M(y,x)-M(0,x)\big),
 \eeq
and $M$ is defined by the following equation:
  \begin{align}
	\begin{aligned}
	K(x,y) &\left(2M(x,y)-M(0,y)\right) = 
		\frac{2x}{3} 
		-2t \by(x+1+\bx)M(x,0)
		+t \by(y+1+\by)M(y,0) \\
		&\qquad
		+t(x-\bx)(y+1+\by)M(0,y)
		-t\left(1 + \by^2 - 2\bx\by \right)%
			M(0,0)
		-t\by 
				M_x(0,0)%
		, 
	\end{aligned}
	\label{eqMMM-king}
\end{align}
where $K(x,y) = 1 - t( x + xy + y + \bx y + \bx + \bx \by + \by + x \by )$.
The specializations $M(x,0)$ and $M(0,y)$ are algebraic each of degree $72$ over $\qs(x,t)$ and $\qs(y,t)$, respectively, and $M(0,0)$ and $M_x(0,0)$ have degree $24$ over $\qs(t)$. 
\end{theorem}

We have moreover a complete algebraic description of all the series
needed to reconstruct $P(x,y)$ and $M(x,y)$ from~\eqref{Psol-king}
and~\eqref{eqMMM-king}, namely the univariate series $M(0,0)$ and
$M_x(0,0)$,
and the bivariate series $M(x,0)$ and $M(0,y)$. 
In particular, both univariate series lie in the extension of $\Q(t)$
(the field of rational functions in $t$) generated in 3 steps as
follows:
first, $\xzero=t+t^2+\LandauO(t^3)$ is the only series in $t$ satisfying
\beq\label{u-def}
  (1-3\xzero)^3 (1+\xzero)t^2+(1 + 18\xzero^2
  -27\xzero^4)t-\xzero=0,
\eeq
then $\av = t+3t^2+\LandauO(t^3)$ is the only series with constant term zero satisfying
\beq\label{v-def}
  (1+3\av-\av^3)\xzero -\av(\av^2+\av+1)=0,
\eeq
and finally
\beq\label{w-def}
  \aw = \sqrt{1 + 4\av - 4\av^3 - 4\av^4} = 1+2t+4t^2+\LandauO(t^3).
\eeq
Schematically,
$	\Q(t) \stackrel{4}{\hookrightarrow} \Q(t,\xzero) \stackrel{3}{\hookrightarrow} \Q(t,\av) \stackrel{2}{\hookrightarrow} \Q(t,\aw). $
Of particular interest is the series $M(0,0)$: by~\eqref{sol-king}, this is
also the series $C_{-1,0}$ that counts by the length walks {in} $\Cc$ ending at
$(-1,0)$. It is algebraic, as conjectured in~\cite{RaschelTrotignon2018Avoiding}, and given by
\beq\label{Cm10-expr}
	M(0,0) =C_{-1,0}= \frac{1}{2t} \left( \frac{\aw(1+2\av)}{1+4\av-2\av^3}-1\right)
	        = t+2t^2+17t^3+80t^4+536t^5+\LandauO(t^6).
\eeq
Due to the lack of space,  the extensions of $\Q(x,t)$ generated by $M(x,0)$ and $M(0,x)$  will only be described in the long version of this paper.	

Once  the series $C(x,y)$ is determined, we can derive detailed
asymptotic results, which refine general results of Denisov and
Wachtel~\cite{DenisovWachtel15} and Mustapha~\cite{Mustapha2019Walks}
(who only obtain the following estimates up to a multiplicative factor).

\begin{corollary}
\label{coro:asy}
The number  $c_{0,0}(n)$ of $n$-step king walks confined to  $\Cc$ and
ending at the origin, and
the number  $c(n)$  of walks of $\Cc$ ending anywhere satisfy for $n \to \infty$:
\begin{align*} 
	c_{0,0}(n) &\sim \left(\frac{2^{29} K}{3^7} \right)^{\! \! 1/3} \,\frac{\Gamma(2/3)}{\pi} \frac{8^n}{n^{5/3}} 
	,\\
	c(n) &\sim \left(\frac{2^{32} K}{3^7}\right)^{\! \! 1/6} \frac{1}{\Gamma(2/3)} \frac{8^n}{n^{1/3}} 
	,
\end{align*}
where 
$K$ is the unique real 
root of $
	101^6 K^3 - 601275603 K^2 + 92811 K - 1$.
%
\end{corollary}

\paragraph*{Outline of the paper} 
We begin in Section~\ref{sec:general} with a general discussion on models of walks with small steps confined to the cone $\Cc$, and on the related functional equations. The main part of the paper, Section~\ref{sec:king}, is devoted to the solution of the king model. We sketch in the final Section~\ref{sec:more} what should be the starting point for the $4$ rightmost models of Figure~\ref{fig:7models}.

\paragraph*{Some definitions and notation}
Let $\GA$ be a commutative ring and $x$ an indeterminate. We denote by
$\GA[x]$ (resp.~$\GA[[x]]$) the ring of polynomials (resp.~\fps) in $x$
with coefficients in $\GA$. If $\GA$ is a field,  then $\GA(x)$ denotes the field of rational functions in $x$, and $\GA((x))$ the
field of Laurent series in $x$, that is, series of the form
$ 
\sum_{n \ge n_0} a_n x^n,
$
with $n_0\in \zs$ and $a_n\in \GA$.  The coefficient of $x^n$ in a   series $F(x)$ is denoted by
$[x^n]F(x)$. 

This notation is generalized to polynomials, fractions,
and series in several indeterminates.
If $F(x,x_1, \ldots, x_d)$ is a series in the $x_i$'s whose
coefficients are Laurent series in $x$, say
\[
F(x,x_1, \ldots, x_d)= \sum_{i_1, \ldots, i_d} x_1^{i_1} \cdots
x_d^{i_d}
\sum_{n \ge n_0(i_1, \ldots, i_d)} a(n, i_1, \ldots, i_d) x^n,
\]
then the \emm non-negative part of $F$ in, $x$ is the
following \fps\ in $x, x_1, \ldots, x_d$:
\[
[x^{\ge 0}]F(x,x_1, \ldots, x_d)= \sum_{i_1, \ldots, i_d} x_1^{i_1} \cdots
x_d^{i_d}
\sum_{n \ge 0} a(n, i_1, \ldots, i_d) x^n.
\]
We define similarly the negative part of $F$, its positive part, and so on.
We denote {with} bars the reciprocals of variables: that is, $\bx=1/x$,
so that $\GA[x,\bx]$ is the ring of Laurent polynomials in $x$ with
coefficients in $\GA$.

If $\GA$ is a field,  a power series $F(x) \in \GA[[x]]$
  is \emm algebraic, (over $\GA(x)$) if it satisfies a
non-trivial polynomial equation $P(x, F(x))=0$ with coefficients in
$\GA$. It is \emm differentially finite, (or \emm D-finite,) if it satisfies a non-trivial linear
differential equation with coefficients  in $\GA(x)$. For
multivariate series, D-finiteness requires the
existence of a differential equation \emm in each variable,.  {We
  refer to~\cite{Li88,lipshitz-df} for general results on D-finite series.}

As mentioned above, we usually omit the dependency in $t$ of our series.
For a series $F(x,y;t) \in \qs[x, \bx, y, \by][[t]]$ and two integers
$i$ and $j$, we denote by $F_{i,j}$ the coefficient of $x^i y^j$ in $F(x,y;t)$. 
This is a series in $\qs[[t]]$.

\section{Enumeration in the three-quarter plane}
\label{sec:general}

We fix a subset $\cS$ of $\{-1,0, 1\}^2\setminus\{(0,0)\}$ and we want
to count walks with steps in $\cS$ that start from the origin $(0,0)$
of $\zs^2$ and remain in the cone $\Cc:=\{(x,y): x\ge 0 \hbox{ or }
y\ge 0\}$. By this, we mean that not only must every vertex of the walk
lie in $\Cc$,
but also every edge: a walk containing a step from
$(-1,0)$ to $(0,-1)$ {(or vice versa)} is not considered as lying in $\Cc$. We often say
for short that our walks \emm avoid the negative quadrant,. The \emm step polynomial, of $\cS$ is defined by
\begin{align*}
  S(x,y)=\sum_{(i,j)\in \cS} x^i y^j%
  = \by H_-(x)+ H_0(x) +y H_+(x)
  = \bx V_-(y) + V_0(y) +x V_+(y),
\end{align*}
for some Laurent polynomials $H_-, H_0, H_+$ and $V_-, V_0, V_+$ (of 
degree {at most} $1$ and valuation {at least} $-1$)
recording horizontal and vertical displacements, respectively.
We denote by $C(x,y;t)\equiv C(x,y)$ the
\gf\ of walks confined to $\Cc$, where the variable $t$ records the
length of the walk, and $x$ and $y$ the coordinates of its endpoints:
\begin{align}\label{C-def}
C(x,y) = \sum_{(i,j) \in \Cc} \sum_{n \geq 0} c_{i,j}(n) x^i y^j t^n
= \sum_{(i,j) \in \Cc} x^i y^j  C_{i,j}(t) .
\end{align}
Here, $c_{i,j}(n)$ is the number of walks of length $n$ that go from
$(0,0)$ to $(i,j)$ and that are confined to $\Cc$.

\subsection{Interesting step sets}
As in the quadrant case~\cite{bomi10}, we can decrease the number of
step sets  that are worth being considered (\emm a priori,, there are $2^8$
of them) thanks to a few simple observations:
\begin{itemize}
\item Since the cone $\Cc$ (as well as the quarter plane $\Qc$) is $x/y$-symmetric, the models defined by $\cS$ and by its mirror image $\overline \cS:=\{(j,i): (i,j) \in \cS\}$ are equivalent; the associated \gfs\ are related by $\overline C(x,y)=C(y,x)$.
\item If all steps of $\cS$ are contained in the right half-plane $\{(x,y): x\ge 0\}$, then \emm all, walks with steps in $\cS$ lie in $\Cc$, and the series $C(x,y)=1/(1-tS(x,y))$ is simply rational. The series $Q(x,y)$ is known to be algebraic in this case~\cite{gessel-factorization}.
\item If all steps of $\cS$ are contained in the left half-plane $\{(x,y): x\le 0\}$, then confining a walk to $\Cc$ is equivalent to confining it to the upper half-plane: the associated \gf\ is then algebraic, and so is $Q(x,y)$.
\item If all steps of $\cS$ lie (weakly) above the first diagonal
  ($x=y$), then confining a walk to~$\Cc$ is again equivalent to confining it to the upper half-plane: the associated \gf\ is then algebraic, and so is $Q(x,y)$.
  \item  If all steps of $\cS$ lie (weakly) above the second diagonal
    ($x+y=0$), then all walks with steps in $\cS$ lie in $\Cc$, and
    $C(x,y)=1/(1-tS(x,y))$ is simply rational. In this case however,
    the series $Q(x,y)$ is not at all
    trivial~\cite{bomi10,Mishna-Rechni}. Such step sets are sometimes
    called \emm singular, in the framework of quadrant walks.
    \item  Finally, if all steps of $\cS$ lie (weakly) below the second diagonal, then a walk confined to~$\Cc$ moves for a while along the second diagonal, and then either stops there or leaves it into the NW or SE quadrant using 
    a South, South-West, or West step. 	It cannot leave the chosen quadrant anymore and behaves therein like a half-plane walk.  
	By polishing this observation, one can prove that $C(x,y)$ is algebraic (while $Q(x,y)=1$).
  \end{itemize}
  Symmetric statements allow us to discard step sets that lie in the upper half-plane $\zs\times \ns$, in the lower half-plane $\zs \times (-\ns)$, or weakly below the $x/y$ diagonal.

 In conclusion, one finds that there are exactly $74$ essentially distinct models of walks avoiding the negative quadrant that are worth studying: the $79$ models considered for quadrant walks (see Tables~$1$--$4$ in~\cite{bomi10}) except the $5$ singular models for which all steps of $\cS$ lie weakly above the diagonal $x+y=0$. 

\subsection{A functional equation}
Constructing walks confined to $\Cc$ step by step gives the following functional equation:
\[
C(x,y)= 1 +tS(x,y) C(x,y) -t\by H_-(x) C_{-,0}(\bx) -t\bx V_-(y)C_{0,-}(\by) -t\bx\by
 C_{0,0} \mathbbm 1_{(-1,-1)\in \cS},
\]
where the series $C_{-,0}(\bx)$ and $C_{0,-}(\by)$ count walks ending on  the horizontal and vertical boundaries of $\Cc$ (but not at $(0,0)$):
\begin{align*}
 \label{Chv}
   C_{-,0}(\bx) &=  \sum_{\substack{i<0\\ n \geq 0}} c_{i,0}(n)x^i t^n \in \bx
   \Q[\bx][[t]],  \\ 
  C_{0,-}(\by) &=  \sum_{\substack{j<0\\ n \geq 0}} c_{0,j}(n)y^j t^n \in \by 
	\Q[\by][[t]].
 \end{align*}
 On the right-hand side of the above functional equation, the term $1$ accounts for the
 empty walk, the next term describes the extension of a walk in $\Cc$
 by one step of $\cS$, and each of the other three terms correspond to
 a ``bad'' move, either starting from the negative $x$-axis, or from
 the negative $y$-axis, or from $(0,0)$. Equivalently,
  \beq\label{eqfunc-gen}
 K(x,y)C(x,y)= 1 -t\by H_-(x) C_{-,0}(\bx) -t\bx V_-(y)C_{0,-}(\by) -t\bx\by
 C_{0,0} \mathbbm 1_{(-1,-1)\in \cS},
 \eeq
 where $K(x,y):=1-tS(x,y)$ is the \emm kernel, of the equation.

 The case of walks confined to the first (non-negative) quadrant
$\Qc$ has been much studied in the past $15$ years. The associated \gf\ $Q(x,y)\equiv Q(x,y;t)\in \qs[x,y][[t]]$ is defined similarly to~\eqref{C-def} and satisfies a similarly looking equation:
 \begin{align*}
 K(x,y)Q(x,y)= 1 -t\by H_-(x) Q_{-,0}(x) -t\bx V_-(y)Q_{0,-}(y) +t\bx\by
   Q_{0,0} \mathbbm 1_{(-1,-1)\in \cS},
 \end{align*}
  where now
\begin{align*}
    Q_{-,0}(x) &=  \sum_{\substack{i\ge0\\ n \geq 0}} q_{i,0}(n)x^i t^n =Q(x,0) \in 
    \Q[x][[t]], \\
		{Q_{0,-}(y)} &=  {\sum_{\substack{j\ge0\\ n \geq 0}} q_{0,j}(n)y^j t^n =Q(0,y) \in 
    \Q[y][[t]].}
\end{align*}

\section{The king walks}
\label{sec:king}
 In this section we focus on the case where the $8$ steps of $\{-1,0,1\}^2\setminus\{(0,0)\}$ are allowed. That is,
 \begin{align*}
	\spol(x,y) = (\bx +1+x)(\by+1+y)-1=x + xy + y + \bx y + \bx + \bx \by + \by + x \by.
\end{align*}
 The functional equation~\eqref{eqfunc-gen} specializes to
\begin{align}
\label{eq:kernelC1}
	K(x,y) C(x,y) &= 1 - t \by(x+1+\bx) C_{-}(\bx) - t \bx(y+1+\by) C_{-}(\by) - t \bx \by C_{0,0},
\end{align}
where we have denoted $ C_{-}(\bx) = C_{-,0}(\bx)= C_{0,-}(\bx)$ (by symmetry).
 Equivalently,
 \begin{align}   \label{eq:kernelC2}
	xy K(x,y) C(x,y) &= xy - t (x^2+x+1) C_{-}(\bx) - t (y^2+y+1) C_{-}(\by) - t C_{0,0}.
\end{align}
The \gf\ $Q(x,y)$ of quadrant walks satisfies
\beq\label{eq:kernelQ1} 
  xy K(x,y) Q(x,y) = xy - t (x^2+x+1) Q(x,0) - t (y^2+y+1) Q(0,y)+ t Q_{0,0}.
\eeq

\subsection{Reduction to an equation with orbit sum zero}
  A key object in the study of walks confined to the first quadrant is
  a certain group of birational transformations that depends on the
  step set. 
	For king walks, it is generated by $(x,y)\mapsto (\bx,
  y)$ and $(x,y)\mapsto (x,\by)$. 
As in~\cite{Bousquet2016}, the similarities between the equations for
$C$ and $Q$, combined with the structure of this group, lead us to
define a new series $A(x,y)$ by
\beq\label{CA-king}
  C(x,y)= A(x,y) + \frac 1 3 \left(  Q(x,y) - \bx^2 Q(\bx,y) - \by^2 Q(x,\by)\right).
\eeq
Then the combination of~\eqref{eq:kernelC2} and~\eqref{eq:kernelQ1} gives
\newcommand*{\myp}{\hspace{-0.15mm}+\hspace{-0.15mm}}
\[ 
xy	K(x,y) A(x,y) = 
		\frac{2xy \myp \bx y \myp x \by}{3} 
		-t(x^2 \myp x \myp 1)A_{-}(\bx)
		-t(y^2 \myp y \myp 1)A_{-}(\by) 
		-tA_{0,0},
\]
and it follows from this equation that  $xyA(x,y)$ has \emph{orbit sum} zero. By this, we mean:
\beq\label{eq:orbitA}
	xy A(x,y) - \bx y A(\bx,y) + \bx \by A(\bx,\by) - x \by A(x,\by) = 0.
\eeq
Theorem~\ref{theo:king} states that $A(x,y)$ is algebraic. 
{In Section~\ref{sec:more} we define 
an analogous series $A$ for all models of
Figure~\ref{fig:7models} which we expect to be systematically
algebraic.}

  The proof of Theorem~\ref{theo:king} starts as in the case of
  the simple and diagonal walks in~\cite{Bousquet2016}. 
	The first objective, achieved in Section~\ref{sec:only}, is to derive an equation that involves a single bivariate series, essentially $A_-(x)$ (and no trivariate series). In principle, the ``generalized quadratic method'' of~\cite{BMJ06} then solves it routinely.
But in practise, the king model turns out to be much {more difficult to solve}
than the other two, and raises 
serious computational difficulties. In what follows, we
focus on the points of the derivation that differ
from~\cite{Bousquet2016}. We have performed all computations with the
computer algebra system {\sc Maple}. The corresponding sessions will be available on the authors' webpages with the long version of the paper.

\subsection{Reduction to a quadrant-like problem}
\label{sec:red_quadrant:king}

We separate in $A(x,y)$ the contributions of the three quadrants, again using the $x/y$-symmetry of the step set: 
\begin{align*}
	A(x,y) = P(x,y) + \bx M(\bx,y) + \by M(\by,x),
\end{align*}
where $P(x,y)$ and $M(x,y)$ lie in $\Q[x,y][[t]]$.
Note that this identity defines $P$ and $M$ uniquely in terms of $A$. 
Replacing $A$ by this expression, and extracting the positive part in $x$ and $y$ from the orbit equation~\eqref{eq:orbitA} relates the series $P$ and $M$ by
\begin{align*}
	xyP(x,y) &= y \left(M(x,y)-M(0,y)\right) + x \left(M(y,x) - M(0,x)\right),
\end{align*}
which is exactly the same as \cite[Eq.~(22)]{Bousquet2016}, and as Eq.~\eqref{Psol-king} in Theorem~\ref{theo:king}.
We then follow the lines of proof of~\cite[Sec.~2.3]{Bousquet2016} to obtain the functional equation~\eqref{eqMMM-king} for~$M$.

\subsection{An equation between \texorpdfstring{$\boldsymbol{M(0,x)$, $M(0,\bx)}$}{M(0,x), M(0,1/x)}, and \texorpdfstring{$\boldsymbol{M(x,0)}$}{M(x,0)}}

Next we will cancel the kernel $K$. As a polynomial in $y$, 
the kernel admits only one root that is a formal power series in $t$:
\begin{align*}
Y(x) = \frac{1 - t(x+\bx) - \sqrt{ (1-t(x+\bx))^2-4t^2(x+1+\bx)^2 }}{ 2t(x+1+\bx)}
	  = (x+1+\bx)t + \LandauO(t^2).
\end{align*}
Note that $Y(x)=Y(\bx)$.
We specialize~\eqref{eqMMM-king} to the pairs $(x,Y(x)), (\bx,Y(x)), (Y(x),x)$, and $(Y(x), \bx)$ (the left-hand side vanishes for each specialization since $K(x,y)=K(y,x)$), and eliminate $M(0,Y)$, $M(Y,0)$, and $M(\bx,0)$  from the four resulting equations. We obtain: 
\begin{align}
\label{eq:3Ms}
\begin{aligned}
	(x+1+\bx)\left(Y(x)-\frac{1}{Y(x)}\right) & \left( xM(0,x) - 2\bx M(0,\bx) \right) 
	+ 3(x+1+\bx)M(x,0) 
	\\ & \qquad
	-\frac{2\bx Y(x)}{t} + 3M_{1,0}
	+ (2Y(x) -x - \bx) M_{0,0} = 0.
\end{aligned}
\end{align}


\subsection{An equation between \texorpdfstring{$\boldsymbol{M(0,x)}$}{M(0,x)} and \texorpdfstring{$\boldsymbol{M(0,\bx)}$}{M(0,1/x)}}

Let us denote the discriminant occurring in $Y(x)$ by
\beq\label{Delta-def}
  \Delta(x) := (1-t(x+\bx))^2-4t^2(x+1+\bx)^2
  = (1-t(3(x+\bx)+2))(1+t(x+\bx+2))
\eeq
and introduce the notation
\begin{align}
\label{eq:RSdef}
\begin{aligned}
	R(x) &:= t^2 M(x,0) = \frac{xt^2}{3} + \left(1+\frac{x^2}{3}\right)t^3 + \LandauO(t^4),  \\
	S(x) &:= txM(0,x) = x(1+x)t^2 + 2x(1+x+x^2)t^3 + \LandauO(t^4).
\end{aligned}
\end{align}
Then~\eqref{eq:3Ms} reads
\begin{align}
\label{eq:3Ms2}
\begin{aligned}
	\sqrt{\Delta(x)} \left( S(x) - 2S(\bx) + \frac{R(0)-t\bx}{t(x+1+\bx)} \right) =
		3(x+1+\bx)&R(x) + 3R'(0)  \\
		  + \frac{1 -t(x+\bx)(x+2+\bx)}{t(x+1+\bx)} &R(0) - \frac{1-t(x+\bx)}{1+x+x^2}.
\end{aligned}
\end{align}
Next, we square this equation and extract the negative part in
$x$. The series $R(x)$ (mostly)  disappears as it involves only non-negative powers of $x$. This gives an expression for the negative part of
$\Delta(x)S(x)S(\bx)$. Using the symmetry of $\Delta(x)$ in $x$ and
$\bx$, we then reconstruct an expression of $\Delta(x)S(x)S(\bx)$ that
does not involve $R(x)$, as in~\cite[Sec.~2.5]{Bousquet2016}. 

During these calculations, we have to extract the negative and
non-negative parts in
series of the form $F(x)/(1+x+\bx)^m$, where $F(x)$ is a series in $t$
with coefficients in $\Q[x,\bx]$. Upon performing a partial fraction
expansion, and separating in $F$ the negative and non-negative parts,
we see that the key question is {how} to {extract and} express 
the non-negative part in series of the form $F({\bx})/(1-\zeta_i x)^m$, {where} $F(x) \in \C[x][[t]]$ and
\[
  \uroot_1 := -\frac{1}{2} + \frac{i\sqrt{3}}{2} \quad \text{ and } \quad
  \uroot_2 := -\frac{1}{2} - \frac{i\sqrt{3}}{2}
\]
are the primitive cubic roots of unity. A simple calculation 
establishes the following lemma.

\begin{lemma}[Non-negative part at pole $\rho$]
\label{lem:pospole}
Let $F(x) \in \C[x][[t]]$ and $\rho \in \C$. Then,
\begin{align*}
		[x^{\geq 0}] \frac{F(\bx)}{1-\rho x} &= \frac{F(\rho)}{1-\rho x}, \\
		[x^{\geq 0}] \frac{F(\bx)}{(1-\rho x)^2} &= \frac{F(\rho)}{(1-\rho x)^2} + \frac{\rho F'(\rho)}{1-\rho x}.
\end{align*}
\end{lemma}
One outcome of the extraction procedure is the following identity:  
\begin{align}
	\label{eq:S1}
	&S(\uroot_1) = S(\uroot_2) = - \frac{R(0) + 3 R'(0)}{1+t} = -t^2 -11t^4 -30t^5 + \LandauO(t^6).
\end{align} 
Using these results, we finally arrive at an equation relating $S(x)$ and $S(\bx)$:
\begin{align}
\label{eq:2Ms}
\begin{aligned}
	\Delta(x) \left( S(x)^2 + S(\bx)^2 - S(x) S(\bx) + \frac{S(x)(xt-R(0)) + \bx S(\bx) ( \bx t - R(0) )}{t(x+1+\bx)} \right) 
	= \\
	(1+t) S(\uroot_1) \left( 2 (x+1+\bx) R(0) - \frac{(1-t(x + \bx))(t(x+\bx)-2R(0))}{t(x+1+\bx)} \right) \\
	+ (1+4t)(x+\bx)R(0) - (t^2 + tR(0) + R(0)^2) (x^2+\bx^2) + \Delta_0,
\end{aligned}
\end{align}
where $\Delta_0$ is the coefficient of $x^0$ in $\Delta(x) S(x)
S(\bx)$.

\subsection{An equation for
  \texorpdfstring{$\boldsymbol{M(0,x)}$}{M(0,x)}
	only}
\label{sec:only}

Equation~\eqref{eq:2Ms} is almost ready for a positive part extraction, except for the mixed term $S(x) S(\bx)$. To eliminate it, we multiply~\eqref{eq:2Ms}
by $S(x) + S(\bx) + \frac{x+\bx-2R(0)/t}{x+1+\bx}$. 
Then we are able to extract the non-negative terms in $x$. 
Hereby we repeatedly apply Lemma~\ref{lem:pospole}. 
Additionally, we use  $R(0)=tS'(0)$ and~\eqref{eq:S1}.
Furthermore, we work with the real and imaginary parts of $\uroot_1 S'(\uroot_1)$ and $\uroot_2 S'(\uroot_2)$. More precisely, we define
\begin{align*}
	(1+t)^2 \uroot_1 S'(\uroot_1) &= \DSA + i\sqrt{3}\DSB, \\
	(1+t)^2 \uroot_2 S'(\uroot_2) &= \DSA - i\sqrt{3}\DSB.
\end{align*}
(Note that $\DSA$ and $\DSB$ here are series in $t$.) 
In the end we get a cubic equation in $S(x)$:
\begin{align}
	\label{eq:Pol1}
	\Pol(S(x), S'(0), S(\zeta_1), \DSA, \DSB, t, x) = 0,
\end{align}
  where the polynomial $\Pol(x_0, x_1, x_2, x_3, x_4 ,t, x)$ is given in Appendix~\ref{app:eqSx}.

\subsection{The generalized quadratic method}
\label{sec:quadratic}
We now use the results of~\cite{BMJ06} to obtain a system of four polynomial equations relating the series  $S'(0), S(\zeta_1), \DSA$, and $\DSB$. Combined with a few initial terms, this system 
characterizes these four series.  Unfortunately, it turned out to be too big for us to solve it completely, be it by bare hand elimination or using Gr\"obner bases: we did obtain a polynomial equation for $S'(0)$ and  $S(\zeta_1)$, but not for the other two series. Instead, we have resorted to a guess-and-check approach,
consisting in \emph{guessing} such equations (of degree $12$ or $24$,
depending on the series), and then \emph{checking}  that they satisfy the
system. This guess-and-check approach is detailed in the next subsection. For the moment, let us explain how the system is obtained.

The approach of~\cite{BMJ06} instructs us to consider the fractional series  $X$ (in $t$), satisfying
\beq\label{dx0}
\Pol_{x_0}(S(X), S'(0), S(\zeta_1), \DSA, \DSB, t, X)
        = 0, 
\eeq
where $\Pol_{x_0}$ stands for the derivative of $\Pol$ with respect to its
first variable. The number and first terms of such series $X$ depend only
on the first  terms of the series $S(x), S'(0), S(\zeta_1), \DSA$, and $\DSB$ (see~\cite[Thm.~2]{BMJ06}). We find that $6$ such series exist:
\begin{align*}
	X_1(t) &= i+2t^2+4t^3+(36-2i)t^4+ \LandauO(t^5), \\
  X_2(t) &= 
           -i+2t^2+4t^3+(36+2i)t^4+ \LandauO(t^5), \\
	X_3(t) &= \sqrt{t}+t+\frac{3}{2}t^{3/2} + 3t^2 + \frac{51}{8}t^{5/2} + 14t^3 + \LandauO(t^{7/2}), \\
	X_4(t) &= -\sqrt{t}+t-\frac{3}{2}t^{3/2} + 3t^2 -
                 \frac{51}{8}t^{5/2} + 14t^3 + \LandauO(t^{7/2}), \\
  	X_5(t) &= i \sqrt{t}-it^{3/2}+2i t^{5/2} + t^3 - 4i t^{7/2} + 2t^4 + \LandauO(t^{9/2}),\\
  X_6(t) &=
           -i \sqrt{t}+it^{3/2}-2i t^{5/2} + t^3 + 4i t^{7/2} + 2t^4 + \LandauO(t^{9/2}).
\end{align*}
Note that the coefficients of $X_1$ and $X_2$ (resp.~$X_5$ and $X_6$)
are conjugates of one another.  As discussed in~\cite{BMJ06}, each of these
series $X$ also satisfies
\beq\label{dx}
\Pol_{x}(S(X), S'(0), S(\zeta_1), \DSA, \DSB, t, X)        = 0, 
\eeq
where $\Pol_x$ is the derivative with respect to the last variable of $\Pol$, and (of course)
\beq\label{original}
\Pol(S(X), S'(0), S(\zeta_1), \DSA, \DSB, t, X)       = 0.
      \eeq
Using this, we can easily identify two of the series $X_i$: indeed,
eliminating $\DSA$ and $\DSB$ between the three equations~\eqref{dx0},~\eqref{dx},
and~\eqref{original} gives a polynomial equation between $S(X), S'(0),  S(\zeta_1),
t$, and $X$, which factors. Remarkably, its simplest non-trivial factor does not
involve $S(X)$, nor $S'(0)$ nor $  S(\zeta_1)$, and reads
\beq\label{X34}
  X^2-t(1+X)^2(1+X^2).
\eeq
By looking at the first terms of the $X_i$'s and the other
factors, one concludes that the above equation holds for $X_3$ and $X_4$, which are thus explicit.

      Let $D(x_1,\ldots,x_4,t,x)$ be the discriminant of $\operatorname{Pol}(x_0,\ldots,x_4,t,x)$ with respect to $x_0$. According to~\cite[Thm.~14]{BMJ06},
      each $X_i$ is a {\emph{double root}} of 
      $D(S'(0),S(\zeta_1), \DSA,\DSB,t,x)$, seen  as a polynomial in $x$.
Hence this polynomial, which involves $4$ unknown series $S'(0),
S(\zeta_1), \DSA,\DSB$, has (at least)
$6$ double roots. This seems  more information than we need! In principle, $4$
double roots should suffice to give $4$ conditions relating  the $4$ unknown
series. However, we shall see that there is some redundancy in the $6$
series $X_i$, which comes from the special form of $D$.

We first observe that $D$ factors as 
\[
  D(S'(0),S(\zeta_1), \DSA,\DSB,t,x)= 27x^2(1+x+x^2)^2 \Delta(x) 
  D_1(S'(0),S(\zeta_1), \DSA,\DSB,t,x),
\]
where $\Delta(x)$ is defined by~\eqref{Delta-def}, and $D_1$ has
degree $24$ in $x$. It is easily checked that none of the $X_i$'s are
roots of the prefactors, so they are double roots of $D_1$. But
we observe that $D_1$ is symmetric in $x$ and $\bx$.  More precisely,
\[
  D_1(S'(0),S(\zeta_1), \DSA,\DSB,t,x)= x^{12} D_2(S'(0),S(\zeta_1),  \DSA,\DSB,t,x+1+\bx), 
\]
for some polynomial $D_2(x_1, \ldots, x_4,t,s)\equiv D_2(s)$  of degree 12 in
$s$. Since each $X_i$ is a double root of $D_1$, each series
$S_i:=X_i+1+1/X_i$, for $1\le i\le 6$,  is a double root of $D_2$. The
series $S_i$, for $2\le i \le 6$,
are easily seen from their first 
terms to be distinct, but the first terms of $S_1$ and $S_2$
suspiciously agree: one suspects (and rightly so), that $X_2=1/X_1$,
and carefully concludes that $D_2$ has (at least) $5$ double roots in
$s$. Moreover, since $X_3$ and $X_4$ satisfy~\eqref{X34}, the
corresponding series $S_3$ and $S_4$ are the roots of $  1+t=tS_i^2$,
that is, $S_{3,4}=\pm \sqrt{1+1/t}$. The other roots start as follows:
\[ 
  S_2= 1+4t^2+8t^3 + \LandauO(t^4),\qquad
  S_{5,6}= \mp \frac i{\sqrt t} + 1+ t^2 \pm i t^{5/2} + \LandauO(t^3).
\]
But this is not the end of the story: indeed, $D_2$ appears to be 
almost symmetric in $s$ and $1/s$. More precisely, we observe that
\[
  D_2(S'(0),S(\zeta_1),  \DSA,\DSB,s)= s^6 D_3\left(S'(0),S(\zeta_1),  \DSA,\DSB,
    ts+\frac{t+1}s\right),
\]
for some polynomial $D_3(S'(0),S(\zeta_1),  \DSA,\DSB,t,z)\equiv D_3(z)$ of degree $6$ in $z$. It follows that each series $Z_i:= tS_i+(1+t)/S_i$, for $2\le i \le 6$,
is a root of $D_3(z)$, and even a double root, unless
$tS_i^2=1+t$, which precisely occurs for $i=3,4$. One finds $Z_{3,4}= \pm 2\sqrt{t(1+t)}$,
\[ 
  Z_2= 1+ 2t -4 t^2 +\LandauO(t^3),\qquad 
  Z_{5,6}= 2t+2t^3 +\LandauO(t^4).
\]
Since $Z_5$ and $Z_6$ seem indistinguishable,  we safely  conclude that $D_3(z)$ has two double roots $Z_2$ and $Z_5$, and a factor $(z^2-4t(1+t))$. 
Writing
\[
  D_3(z)= \sum_{i=0}^6 d_i z^i,
\]
these properties imply, by matching the three monomials of highest degree, that
\[
  D_3(z)= {\frac { \left({z}^{2}-4\,t(1+t) \right) \left( 8\,{z}^{2}{d_{{6}}}^{2}+4\,zd_{{5}}d_{{
  6}}+16\,{t}^{2}{d_{{6}}}^{2}+16\,t{d_{{6}}}^{2}+4\,d_{{4}}d_{{6}}-{d_{{5}}}^{2} \right) ^{2}}{64\,{d_{{6}}}^{3}}}.
\]
Extracting the coefficients of $z^0, \ldots, z^3$ gives $4$ polynomial
relations between  the coefficients~$d_i$, resulting in $4$ polynomial
relations between the $4$ series $S'(0),S(\zeta_1),  \DSA,\DSB$. One easily
checks that this system, combined with the first terms of these
series, defines them uniquely.

As explained at the beginning of this subsection, we have at the moment only been able to derive from this system polynomial
equations (of degree $24$) for  $S'(0)$ and $S(\zeta_1)$. For the other
two, we had to resort to a guess-and-check approach, which we now describe.

%
\subsection{Guess-and-check}
\label{sec:guess}
{\bf Guessing.}
Returning to the functional equation~\eqref{eq:kernelC1} it is easy to
extract a simple recurrence for the polynomials $c_n(x,y)$
that count walks of length $n$ by the position of their endpoint. 
We implemented this recurrence in the programming language $C$ using {modular} arithmetic and the Chinese remainder theorem to compute the explicit values of this sequence up to $n=2000$. 
Then we were able to guess polynomial 
equations satisfied by $S'(0)$, $S(\zeta_1)$, $\DSA$, and $\DSB$ {using the \texttt{gfun} package in {\sc Maple}~\cite{gfun}}. Of course, those obtained for $S'(0)$ and $S(\zeta_1)$ coincide with those that we derived from the system of the previous subsection.
Details on the corresponding equations are shown below.


\begin{center}
  \begin{tabular}{|c|c|c|c|}
    \hline
    Generating function & Degree in $GF$ & Degree in $t$  & Number of terms \\
    \hline	
		$S'(0)$ &  $24$ & $12$ & $323$ \\
		$S(\zeta_1)$ &  $24$ & $32$ & $823$ \\
		$\DSA$ &  $12$ & $26$ & $229$ \\
		$\DSB$ &  $24$ & $60$ & $477$ \\
    \hline
  \end{tabular}
\end{center}

Checking  that the guessed series satisfy the system turns out to be much easier once the algebraic structure of these series is elucidated, which we do below\footnote{For this section, we have greatly benefited from the help of Mark van Hoeij (\url{https://www.math.fsu.edu/~hoeij/}), who  explained us how to find subextensions of $\Q(t,\DSA)$, and ``simple'' series in these extensions.}. We have not tried a direct check.

\newcommand{\Pmab}{P}


\bigskip \noindent
{\bf The algebraic structure of 
  $\boldsymbol{S'(0), S(\zeta_1), \DSA}$, and $\boldsymbol{\DSB}$.}
%
We begin with the simplest series,~$\DSA$, of (conjectured) degree
$12$. Let $\Pmab(F,t)$ be its guessed monic minimal polynomial.
Using the \texttt{Subfields} command of
{\sc Maple} for several fixed values of $t$, one conjectures that the extension
$\Q(t,\DSA)$ possesses a subfield $Q(t,\xzero)$ of degree $4$ over $\Q(t)$.
{\sc Maple} gives a possible generator $\xzero$ for fixed values of $t$, but how
can we choose $u$ for a generic $t$? Indeed, the value of $\xzero$
given by {\sc Maple} for fixed $t$ has 
no reason to be 
canonical. But the factorisation of $\Pmab(F,t)$ over $\Q(t,\xzero)$, of
the form $P_3(F) P_9(F)$ (with $P_i$ of degree $i$), with coefficients in $\Q(t,\xzero)$, {\emph{is}} canonical. Hence we will compute this factorisation, first
for fixed values of $t$. We proceed as follows: we factor $\Pmab(F,t)$
over $\Q(t,\DSA)$, and find,  for fixed $t=3,\ldots,50$, that
\begin{align*}
	\Pmab(F,t) =  (F-\DSA) P_2(F,\DSA) P_9(F,\DSA),
\end{align*}
where $P_2$ (resp.~$P_9$) is a monic polynomial of degree $2$ (resp.~$9$)  in
$F$. Hence the cubic factor $P_3(F)= F^3+p_2 F^2+p_1 F+p_0$ must be
$(F-\DSA) P_2(F,\DSA)$, and we have just found its
coefficients $p_i$ in terms of $\DSA$ (for $t$ fixed). We now compute the minimal polynomial over $\Q$ of each $p_i$ using a resultant or the \texttt{evala/Norm} command in {\sc Maple}. If 
the above factorization persists for all $t$,
as we expect, each $p_i$ should have a minimal polynomial over $\Q(t)$
of degree (at most) $4$. Having computed this polynomial for
sufficiently many values of~$t$, we reconstruct its generic form by
rational reconstruction. We find that all $p_i$ generate the same
extension of degree $4$ of $\Q(t)$, and we can take any of them as a
first candidate for the generator $\xzero$.
We may  simplify this generator further to end with the
choice~\eqref{u-def}. Then we factor $\Pmab(F,t)$ over $\Q(t,
\xzero)$, and check that our guess was correct: 
the series $\DSA$ is indeed cubic over $\Q(t,\xzero)$. Moreover, it can be written  rationally  in terms of  $t$
and the series $\av$ given by~\eqref{v-def}.

 Finally, we factor the guessed minimal polynomials of $S'(0)$,
$S(\zeta_1)$, and $\DSB$ over $\Q(t,\av)$, and find that these three series all belong to the same quadratic extension of $\Q(t,\av)$, generated by the series
$w$ given by~\eqref{w-def}. In particular,
\[
  S'(0)= \frac 1 2 \left(  \frac{\aw(1+2\av)}{1+4\av-2\av^3}-1\right),
\]
which coincides with~\eqref{Cm10-expr}, given the Definition~\eqref{eq:RSdef} of $S(x)$.

Now that we have guessed rational expressions of $S'(0)$, $S(\zeta_1)$, $\DSA$, and  $\DSB$ in terms of $t$, $\av$, and $\aw$, the $4$ equations obtained in  Section~\ref{sec:quadratic}       
are readily checked to hold, using  the minimal polynomials of $\av$ and $\aw$.

\subsection{Back to \texorpdfstring{$\boldsymbol{S(x)}$}{S(x)} and \texorpdfstring{$\boldsymbol{R(x)}$}{R(x)}}

For $S(x)$ we start with Equation~\eqref{eq:Pol1}, with all one-variable series replaced by their expressions in terms of $t$, $v$, and $w$. 
We eliminate $\aw$ and $\av$ using resultants to arrive at an equation of degree~$72$ over $\Q(t,x)$ for $S(x)= txM(0,x)$.

We can simplify~\eqref{eq:Pol1}  by      
working with the depressed equation, i.e., removing the quadratic term by a suitable change of variable. Indeed, defining $T(x)$ by
\begin{align*}
		S(x) = T(x) + \frac{3xS'(0)-2x^2-1}{3(x^2+x+1)},
\end{align*} 
we find that $T(x)$ satisfies a cubic equation with no quadratic term, involving $t$ and $\av$ but not $\aw$. That is, $T(x)$ has degree $36$ over $\Q(t,x)$, instead of $72$ for $S(x)$.

Introducing $T(x)$ also helps understanding the algebraic structure of  $R(x)$. Returning to~\eqref{eq:3Ms2}, we  recall that $R(0)=tS'(0)$ and use~\eqref{eq:S1} to express $R'(0)$ in terms of $t, v$, and $w$.
The left-hand side simply reads $\sqrt {\Delta(x)}(T(x) - 2T(\bx))$, and is found to be an element of $w \Q(t,x, T(x))$. In the end, $R(x)$ has degree $72$ and belongs to the same extension of $\Q(t,x)$  as $S(x)$.
{This ends the proof of our main result, Theorem~\ref{theo:king}.}

%
%

\section{More models}
\label{sec:more}
For each of the $7$ step sets $\Sc$ of Figure~\ref{fig:7models}, we are able to define a series $A(x,y)$ that
\begin{itemize}
  \item satisfies the same equation as $C(x,y)$ (see~\eqref{eqfunc-gen}), but with a different constant term,
  \item satisfies an \emph{orbit sum} identity similar to~\eqref{eq:orbitA}.
\end{itemize}
Explaining where this series comes from would require {us} to introduce the group associated to a step set. For the sake of conciseness, we simply define $A(x,y)$ without further justification.

For the first four step sets $\Sc$  of Figure~\ref{fig:7models},
the series $A(x,y)$ is defined by~\eqref{CA-king} (with $Q(x,y)$ counting quadrant walks with steps in $\Sc$) as we have seen. 
For the next two step sets, 
\[ 
  C(x,y)=A(x,y)+ \frac 1 5\left( Q(x,y)-\bx^2y Q( \bx y, y) + \bx^3
   Q(\bx y, \bx) + \by^3 Q(\by, x\by) -x\by^2Q(x,x\by)\right).
\] 
 Finally, for the seventh one, 
\begin{multline*}  
C(x,y)=A(x,y)+ \frac 1 7\left( Q(x,y)
   - \bx^2 y Q(\bx y,y)+  \bx^4 y Q(\bx y,\bx^2   y) \right.
 \\ \left.
   -\bx^4 Q(\bx,\bx^2 y )  -\by^3 Q  (x\by ,\by)+  x^2 \by^3 Q (x\by ,x^2 \by) -x^2 \by^2 Q  (x ,x^2 \by)\right)   .
\end{multline*}
In all cases, the series $A(x,y)$ satisfies the following variant of~\eqref{eqfunc-gen}:
\[
 K(x,y)A(x,y)= P_0(x,y) -t\by H_-(x) A_{-,0}(\bx) -t\bx V_-(y)A_{0,-}(\by) -t\bx\by
 A_{0,0} \mathbbm 1_{(-1,-1)\in \cS},
\]
where $K(x,y)=1-tS(x,y)$ as before, and $P_0(x,y)$ is a Laurent polynomial. This equation is easily obtained by combining the equations for $C(x,y)$ and $Q(x,y)$.

Finally, the vanishing orbit sum, which is~\eqref{eq:orbitA} for the first four  models, reads
\[ 
  xyA(x,y)- \bx y^2 A(\bx y,y)+  \bx^2 y A(\bx y,\bx)- \bx \by A
  (\by,\bx) + x\by^2 A (\by,x\by)- x^2\by A(x,x\by)=0
\] 
for the next two, and
 \begin{multline*}
 xyA(x,y)- \bx y^2 A(\bx y,y)+  \bx^3 y^2 A(\bx y,\bx^2   y)
  -\bx^3 y A(\bx,\bx^2 y )
 \\
 + \bx\by A(\bx,\by) -x\by^2 A  (x\by ,\by)
+ x^3 \by^2 A (x\by ,x^2 \by) -x^3 \by A  (x ,x^2 \by) =0
  \end{multline*} 
for the last one.
We conjecture that the series $A(x,y)$ is systematically algebraic (this is now proved for the first three models). To support this conjecture, we have tried to guess {(using the \texttt{gfun} package~\cite{gfun} in {\sc Maple})}, for the~$4$ models for which it is still open, a polynomial equation for the series $A_{-1,0}$, which, in all cases, coincides with the \gf\ $C_{-1,0}$ of walks ending at $(-1,0)$ (for the second model we consider $A_{-2,0}$ instead, since $A_{-1,0}=0$ due to the periodicity of the model).  This series has degree $4$ (resp.~$8$, $24$)  in the three solved cases. We could not guess anything for the $4$th model {(using the counting sequence for such walks up to length $n=4000$)}, but we discovered equations of degree $24$ for each of the next three.

      We believe that it would be worth exploring if the guiding principles of the present paper apply to these $4$ other models. In all cases, we expect to face a \emph{system} of quadrant-like equations rather than a single one. We plan to investigate at least some of these models.
  
      To conclude, we recall that the $4$ small step models that are algebraic for the quadrant problem are conjectured to be algebraic for the three-quadrant cone as well~\cite[Fig.~5]{Bousquet2016}.  In this case, the series $A(x,y)$ simply coincides with $C(x,y)$, as the orbit sum of $xyC(x,y)$ vanishes.


\appendix

\section{Final polynomial equation for \texorpdfstring{$\boldsymbol{S(x)}$}{S(x)} in the king model}
\label{app:eqSx}

The polynomial $\Pol$ involved in the cubic Equation~\eqref{eq:Pol1}  defining  $S(x)$ is: 
\begin{align*}
	&\operatorname{Pol}(x_0, x_1, x_2, x_3, x_4, t, x) = \\
	&-3(x^2+x+1)^2(x^2t+2xt+x+t)(3x^2t+2xt-x+3t) \, \text{\textcolor{blue}{\boldmath $x_0^3$}}\\
	&\quad +3(x^2+x+1)(x^2t+2xt+x+t)(3x^2t+2xt-x+3t)(3x_1x-2x^2-1) \, \text{\textcolor{blue}{\boldmath $x_0^2$}}\\
	&\quad + \left[ 3x^2(x^2+x+1)^2(2x_4x_1+x_4-x_3)-3x^2(t+1)^2(x^2+x+1)^2x_2^2 \right.\\
	&\quad + 6x(t+1)(x^2+x+1)(x^4t+2x^2t+x^2+t)x_1x_2 \\
	&\quad +3x(t+1)(x^2+x+1)(x^4t-x^3t-x^3+x^2t-xt-x+t)x_2 
	-3\left(x^8t^2+2x^7t^2 \right.\\
	&\quad \left.+10x^6t^2+20x^5t^2+4x^5t+25x^4t^2+20x^3t^2-2x^4+4x^3t+10x^2t^2+2xt^2+t^2\right)x_1^2 \\
	&\quad -3\left(x^8t^2-11x^7t^2-x^7t-32x^6t^2-9x^6t-53x^5t^2-6x^5t-55x^4t^2+3x^5-15x^4t \right.\\
	&\quad \left. -39x^3t^2-6x^3t-16x^2t^2+x^3-5x^2t-5xt^2-xt+t^2\right)x_1-12x^8t^2-30x^7t^2-6x^7t \\
	&\quad   -51x^6t^2-60x^5t^2+3x^6-12x^5t-54x^4t^2-36x^3t^2+3x^4-6x^3t -21x^2t^2-6xt^2\\
	&\quad \left. -3t^2 \right]\text{\textcolor{blue}{\boldmath $x_0$}} + x^2(x^2+x+1)\left[(2x_3x^2-6x_4x-2x_3)x_1^2  -(x^2+2)x_3+3x_4x^2  \right. \\
	&\quad \left. +(2x-1)(3x_4x+x_3 (x+2))x_1\right] +3x^3(t+1)^2(x^2+x+1)(x_1-x)x_2^2 \\
	&\quad -3x^2(t+1)x_2 (x_1-x)((2(x^2+t(x^2+1)^2))x_1+t(x^4+x^2+1) -(t+1)x(x^2+1))\\
	&\quad +3xt(x^2+x+1)^2(x_1-x)(t(x^2-x+1)x_1^2+(x^2t-5xt-x+t)x_1+t(x^2-x+1)).
\end{align*}



\bibliography{Bibliography}

\begin{thebibliography}{10}

\bibitem{Bousquet2016}
M.~Bousquet-M\'{e}lou.
\newblock Square lattice walks avoiding a quadrant.
\newblock {\em J. Combin. Theory Ser. A}, 144:37--79, 2016.
\newblock \href{http://arxiv.org/abs/1511.02111}{arXiv:1511.02111}.
\newblock \href {https://doi.org/10.1016/j.jcta.2016.06.010}
  {\path{doi:10.1016/j.jcta.2016.06.010}}.

\bibitem{BMJ06}
M.~Bousquet-M\'elou and A.~Jehanne.
\newblock Polynomial equations with one catalytic variable, algebraic series
  and map enumeration.
\newblock {\em J. Combin. Theory Ser. B}, 96:623--672, 2006.
\newblock \href{http://arxiv.org/abs/math/0504018}{arXiv:math/0504018}.
\newblock \href {https://doi.org/10.1016/j.jctb.2005.12.003}
  {\path{doi:10.1016/j.jctb.2005.12.003}}.

\bibitem{bomi10}
M.~Bousquet-M\'{e}lou and M.~Mishna.
\newblock Walks with small steps in the quarter plane.
\newblock In {\em Algorithmic probability and combinatorics}, volume 520 of
  {\em Contemp. Math.}, pages 1--39. Amer. Math. Soc., Providence, RI, 2010.
\newblock \href{http://arxiv.org/abs/0810.4387}{arXiv:0810.4387}.

\bibitem{DenisovWachtel15}
D.~Denisov and V.~Wachtel.
\newblock Random walks in cones.
\newblock {\em Ann. Probab.}, 43(3):992--1044, 2015.
\newblock \href{http://arxiv.org/abs/1110.1254}{arXiv:1110.1254}.
\newblock \href {https://doi.org/10.1214/13-AOP867}
  {\path{doi:10.1214/13-AOP867}}.

\bibitem{fayolle-livre}
G.~Fayolle, R.~Iasnogorodski, and V.~Malyshev.
\newblock {\em Random walks in the quarter-plane: Algebraic methods, boundary
  value problems and applications}, volume~40 of {\em Applications of
  Mathematics}.
\newblock Springer-Verlag, Berlin, 1999.
\newblock \href{https://doi.org/10.1007/978-3-319-50930-3}{[doi]}.

\bibitem{gessel-factorization}
I.~Gessel.
\newblock A factorization for formal {L}aurent series and lattice path
  enumeration.
\newblock {\em J. Combin. Theory Ser. A}, 28(3):321--337, 1980.
\newblock \href{https://doi.org/10.1016/0097-3165(80)90074-6}{[doi]}.
\newblock \href {https://doi.org/10.1016/0097-3165(80)90074-6}
  {\path{doi:10.1016/0097-3165(80)90074-6}}.

\bibitem{gessel-zeilberger}
I.~M. Gessel and D.~Zeilberger.
\newblock Random walk in a {W}eyl chamber.
\newblock {\em Proc. Amer. Math. Soc.}, 115(1):27--31, 1992.
\newblock \href{https://doi.org/10.1090/S0002-9939-1992-1092920-8}{[doi]}.
\newblock \href {https://doi.org/10.1090/S0002-9939-1992-1092920-8}
  {\path{doi:10.1090/S0002-9939-1992-1092920-8}}.

\bibitem{Li88}
L.~Lipshitz.
\newblock The diagonal of a {$D$}-finite power series is {$D$}-finite.
\newblock {\em J. Algebra}, 113(2):373--378, 1988.
\newblock \href{https://doi.org/10.1016/0021-8693(88)90166-4}{[doi]}.
\newblock \href {https://doi.org/10.1016/0021-8693(88)90166-4}
  {\path{doi:10.1016/0021-8693(88)90166-4}}.

\bibitem{lipshitz-df}
L.~Lipshitz.
\newblock D-finite power series.
\newblock {\em J. Algebra}, 122:353--373, 1989.
\newblock \href{https://doi.org/10.1016/0021-8693(89)90222-6}{[doi]}.
\newblock \href {https://doi.org/10.1016/0021-8693(89)90222-6}
  {\path{doi:10.1016/0021-8693(89)90222-6}}.

\bibitem{Mishna-Rechni}
M.~Mishna and A.~Rechnitzer.
\newblock Two non-holonomic lattice walks in the quarter plane.
\newblock {\em Theoret. Comput. Sci.}, 410(38-40):3616--3630, 2009.
\newblock \href{http://arxiv.org/abs/math/0701800}{arXiv:math/0701800}.
\newblock \href {https://doi.org/10.1016/j.tcs.2009.04.008}
  {\path{doi:10.1016/j.tcs.2009.04.008}}.

\bibitem{Mustapha2019Walks}
S.~Mustapha.
\newblock Non-{D}-finite walks in a three-quadrant cone.
\newblock {\em Ann. Comb.}, 23(1):143--158, 2019.
\newblock \href{https://doi.org/10.1007/s00026-019-00413-2}{[doi]}.
\newblock \href {https://doi.org/10.1007/s00026-019-00413-2}
  {\path{doi:10.1007/s00026-019-00413-2}}.

\bibitem{raschel-unified}
K.~Raschel.
\newblock Counting walks in a quadrant: a unified approach via boundary value
  problems.
\newblock {\em J. Eur. Math. Soc. (JEMS)}, 14(3):749--777, 2012.
\newblock \href{http://arxiv.org/abs/1003.1362}{arXiv:1003.1362}.
\newblock \href {https://doi.org/10.4171/JEMS/317}
  {\path{doi:10.4171/JEMS/317}}.

\bibitem{RaschelTrotignon2018Avoiding}
K.~Raschel and A.~Trotignon.
\newblock On walks avoiding a quadrant.
\newblock {\em Electron. J. Combin.}, 26(3):\href{
  https://doi.org/10.37236/8019}{Paper 3.31}, 34 pp., 2019.
\newblock \href{https://arxiv.org/abs/1807.08610}{arXiv:1807.08610}.

\bibitem{gfun}
B.~Salvy and P.~Zimmermann.
\newblock Gfun: a {M}aple package for the manipulation of generating and
  holonomic functions in one variable.
\newblock {\em ACM Transactions on Mathematical Software}, 20(2):163--177,
  1994.
\newblock \href{https://doi.org/10.1145/178365.178368}{[doi]}.
\newblock \href {https://doi.org/10.1145/178365.178368}
  {\path{doi:10.1145/178365.178368}}.

\end{thebibliography}

\end{document}